%% file: main.tex
\documentclass{informs3}

\OneAndAHalfSpacedXI

\usepackage{natbib}
 \bibpunct[, ]{(}{)}{,}{a}{}{,}%
 \def\bibfont{\small}%
\TheoremsNumberedThrough     
\ECRepeatTheorems
\EquationsNumberedThrough    

\input{packages.tex}

\input{macros.tex}

\begin{document}

\def\COPYRIGHTHOLDER{INFORMS}%
\def\COPYRIGHTYEAR{2017}%
\def\DOI{\fontsize{7.5}{9.5}\selectfont\sf\bfseries\noindent https://doi.org/10.1287/opre.2017.1714\CQ{Word count = 9740}}

\RUNAUTHOR{Lan and Berkhout}

\RUNTITLE{PyJobShop: Solving scheduling problems with constraint programming in Python}

\TITLE{PyJobShop: Solving scheduling problems with constraint programming in Python}

\ARTICLEAUTHORS{
	\AUTHOR{Leon Lan, Joost Berkhout}
	\AFF{Department of Mathematics, Vrije Universiteit Amsterdam, \EMAIL{\{l.lan,joost.berkhout\}@vu.nl}}}

    \input{sections/abstract}



\maketitle

\input{sections/introduction}
\input{sections/problem}

\input{sections/variants}
\input{sections/software}

\input{sections/experiments}
\input{sections/conclusion}

\ACKNOWLEDGMENT{This work was supported by TKI Dinalog, Topsector Logistics, and the Dutch Ministry of Economic Affairs and Climate Policy.
}	

\bibliographystyle{informs2014}
\bibliography{references.bib}

\begin{APPENDICES}
  \input{appendices/experiments}
\end{APPENDICES}

\end{document}

%% file: packages.tex
\usepackage{anyfontsize}

\usepackage{hyperref}
\usepackage{multicol}
\usepackage{enumitem} 
\usepackage{bm} 
\usepackage{bbm} 
\usepackage{amsmath,amssymb}
\usepackage{mathrsfs} 

\usepackage{listings}
\usepackage{pythonhighlight}
\usepackage{tabularx}
\usepackage{booktabs}
\usepackage{multirow}
\usepackage{siunitx}

%% file: macros.tex
\NewDocumentCommand{\StartAtStart}{ o o }{%
  \texttt{StartAtStart}%
  \IfValueT{#1}{(#1, #2)}%
}

\NewDocumentCommand{\StartAtEnd}{ o o }{%
  \texttt{StartAtEnd}%
  \IfValueT{#1}{(#1, #2)}%
}

\NewDocumentCommand{\StartBeforeStart}{ o o }{%
  \texttt{StartBeforeStart}%
  \IfValueT{#1}{(#1, #2)}%
}

\NewDocumentCommand{\StartBeforeEnd}{ o o }{%
  \texttt{StartBeforeEnd}%
  \IfValueT{#1}{(#1, #2)}%
}

\NewDocumentCommand{\EndAtStart}{ o o }{%
  \texttt{EndAtStart}%
  \IfValueT{#1}{(#1, #2)}%
}

\NewDocumentCommand{\EndAtEnd}{ o o }{%
  \texttt{EndAtEnd}%
  \IfValueT{#1}{(#1, #2)}%
}

\NewDocumentCommand{\EndBeforeStart}{ o o }{%
  \texttt{EndBeforeStart}%
  \IfValueT{#1}{(#1, #2)}%
}

\NewDocumentCommand{\EndBeforeEnd}{ o o }{%
  \texttt{EndBeforeEnd}%
  \IfValueT{#1}{(#1, #2)}%
}

\NewDocumentCommand{\SameUnit}{ o o }{%
  \texttt{SameUnit}%
  \IfValueT{#1}{(#1, #2)}%
}

\NewDocumentCommand{\DifferentUnits}{ o o }{%
  \texttt{DifferentUnits}%
  \IfValueT{#1}{(#1, #2)}%
}

\NewDocumentCommand{\Previous}{ o o }{%
  \texttt{Previous}%
  \IfValueT{#1}{(#1, #2)}%
}

\NewDocumentCommand{\Accessibility}{ o o }{%
  \texttt{Accessibility}%
  \IfValueT{#1}{(#1, #2)}%
}

\NewDocumentCommand{\StartOf}{o}{%
  \texttt{StartOf}%
  \IfValueT{#1}{(#1)}%
}

\NewDocumentCommand{\EndOf}{o}{%
  \texttt{EndOf}%
  \IfValueT{#1}{(#1)}%
}

\NewDocumentCommand{\PresenceOf}{o}{%
  \texttt{PresenceOf}%
  \IfValueT{#1}{(#1)}%
}

\newcommand{\NoOverlap}{\texttt{NoOverlap}}
\newcommand{\Cumulative}{\texttt{Cumulative}}

\newcommand{\start}{\text{start}}
\newcommand{\supend}{\text{end}}
\newcommand{\duration}{\text{duration}}
\newcommand{\present}{\text{present}}

%% file: sections/abstract.tex
\ABSTRACT{%
    This paper presents PyJobShop, an open-source Python library for solving scheduling problems with constraint programming.
    PyJobShop provides an easy-to-use modeling interface that supports a wide variety of scheduling problems, including well-known variants such as the flexible job shop problem and the resource-constrained project scheduling problem.
    PyJobShop integrates two state-of-the-art constraint programming solvers: Google's OR-Tools CP-SAT and IBM ILOG's CP Optimizer.
    We leverage PyJobShop to conduct large-scale numerical experiments on more than 9,000 benchmark instances from the machine scheduling and project scheduling literature, comparing the performance of OR-Tools and CP Optimizer.
    While CP Optimizer performs better on permutation scheduling and large-scale problems, OR-Tools is highly competitive on job shop scheduling and project scheduling problems--while also being fully open-source.
    By providing an accessible and tested implementation of constraint programming for scheduling, we hope that PyJobShop will enable researchers and practitioners to use constraint programming for real-world scheduling problems.
}

%% file: sections/introduction.tex
\section{Introduction}
Scheduling is a key process in both manufacturing and service sectors, as it involves the efficient allocation of resources to tasks over time. 
Challenges in scheduling span many diverse domains--from semiconductor manufacturing to construction project planning--making it one of the most intensively researched topics in operations research \citep{potts_fifty_2009}. 
Over the past several decades, researchers have developed a wide range of models, evolving from simple single-machine settings to more complex job shop environments, with solution techniques that vary from straightforward dispatching rules to advanced metaheuristic algorithms.

Constraint programming (CP), a technique rooted in the field of artificial intelligence, has emerged as a promising technique to address scheduling problems, outperforming the traditional mixed-integer linear programming (MILP) approach~\citep{ku_mixed_2016, naderi_mixed-integer_2023}.
Its interval-based modeling language results in formulations that are more intuitive and compact than their MILP counterparts, while constraint propagation combined with search techniques make CP highly effective at obtaining high-quality solutions for scheduling \citep{laborie_ibm_2018}.
This effectiveness has led to the widespread adoption of CP in solving scheduling problems over the past decade \citep{naderi_mixed-integer_2023}.

Motivated by this increasing interest in CP, we present PyJobShop, an open-source Python package to solve scheduling problems with CP.
PyJobShop provides an easy-to-use modeling interface that allows users to solve a large variety of machine and project scheduling problems without having to understand the specific CP implementation details.
PyJobShop implements a general scheduling model and integrates two state-of-the-art CP solvers: Google's OR-Tools CP-SAT and IBM ILOG's CP Optimizer.
The package features comprehensive documentation, extensive testing, and extensible CP models, and our code is open-source under the permissive MIT License.
The contributions of our work can be further categorized as follows.

\begin{description}
\item[General-purpose scheduling model implemented with CP.] 
Inspired by the review on the flexible job shop problem by \citet{dauzere-peres_flexible_2024}, we introduce a general-purpose scheduling model that forms the basis of PyJobShop. 
This approach provides a single interface to address a wide variety of scheduling variants, including machine scheduling problems and also variants from the project scheduling literature.
PyJobShop's scheduling model is implemented using Google's OR-Tools CP-SAT \citep{perron_or-tools_2024} and IBM ILOG's CP Optimizer \citep{laborie_ibm_2018}. 
In line with academic conventions, we refer to OR-Tools CP-SAT as OR-Tools and IBM ILOG's CP Optimizer as CP Optimizer.
OR-Tools is an open-source solver that integrates CP and SAT solving techniques through lazy clause generation \citep{stuckey_lazy_2010} and has demonstrated state-of-the-art performance by consistently achieving first place in the MiniZinc CP challenge since 2018 \citep{minizinc_minizinc_2025}.
In contrast, CP Optimizer is a commercial solver that stems from ILOG’s long history in CP \citep{baptiste_constraint-based_2001}.
It is widely adopted within the scheduling community and is recognized for its strong performance, having successfully closed many long-standing open instances \citep{laborie_ibm_2018}.

\item[Performance comparison between OR-Tools and CP Optimizer.]
We use PyJobShop to conduct a large-scale numerical experiment to compare the performance between OR-Tools and CP Optimizer using more than 9,000 benchmark instances from the machine scheduling and project scheduling literature.
Our results show that OR-tools is highly competitive with CP Optimizer —delivering strong results on job shop and project scheduling problems while also achieving superior lower bounds.
CP Optimizer, on the other hand, demonstrates superior scalability for permutation scheduling problems and large-scale problems.
In addition, both OR-Tools and CP Optimizer find many new best-known solutions on benchmark instances from the project scheduling literature, further showcasing the strengths of CP in addressing these types of scheduling problems.

\item[Open-source, tested and documented.]
Our final contribution is an open-source, tested, and documented software package that addresses a significant gap in the scheduling literature, where many implementations remain unpublished and are difficult to reproduce.
We are aware of two related scheduling software projects: SSP-3M \citep{marquez_open-source_2024}, an open-source framework for shop scheduling problems focused on designing heuristics, and the Job Shop Scheduling Benchmark \citep{reijnen_job_2023}, a benchmark library focused on shop scheduling problems with special emphasis on reinforcement learning methods.
PyJobShop distinguishes itself through its comprehensive documentation, extensive test coverage, modern CP solver integration, and support for multiple scheduling variants through a single interface.
Our focus on software makes PyJobShop particularly suitable for both research extensions and practical applications, and all code, data, and results are open-sourced to promote reproducibility and future development.
\end{description}

The outline of the remaining paper is as follows.
In Section~\ref{sec:problem}, we describe the scheduling model of PyJobShop and present its CP formulation.
Section~\ref{sec:variants} gives an overview of all problem variants that PyJobShop supports.
The PyJobShop software package is presented and demonstrated in Section~\ref{sec:software}.
Section~\ref{sec:experiments} describes and presents the large-scale numerical experiments.
The paper is concluded in Section~\ref{sec:conclusion}.

%% file: sections/problem.tex
\section{Problem description}
\label{sec:problem}
In this section, we describe PyJobShop's scheduling model.
In Section~\ref{subsec:problem-notation}, we introduce the notation and preliminaries and informally describe the problem.
In Section~\ref{subsec:model}, we present the CP formulation.
In the following, we assume that all numerical values are integral because CP solvers generally do not support continuous values.

\subsection{Problem notation and preliminaries}
\label{subsec:problem-notation}
Let $J$ be the set of jobs, $R$ the set of resources, $T$ the set of tasks, and $M$ the set of modes.
A job $j \in J$ represents a collection of tasks whose completion influences the objective.
Each job $j \in J$ has a set of related tasks $T_{j} \subseteq T$, a release date $r_j \geq 0$ when the job becomes available, a deadline $\overline{d_{j}} \geq 0$ by which the job must be completed, and an optional due date $d_j \geq 0$ by which the job is ideally completed, and a weight $w_j \geq 0$ reflecting its priority.

A resource $r \in R$ is used to process tasks, and the set of resources $R$ is partitioned into three disjoint sets $R =  R^{\text{machine}} \cup R^{\text{renewable}} \cup R^{\text{non-renewable}}$.
A machine $r \in R^{\text{machine}}$ is a resource that can process only one task at a time and can handle sequencing constraints.
A renewable resource $r \in R^\text{renewable}$ is a resource that has capacity $Q_r \ge 0$ at each time period.
In contrast, a non-renewable resource $r \in R^\text{\text{non-renewable}}$ is a resource that can process at most $Q_r \ge 0$ demand in total over the entire time horizon.

A task $t \in T$ is the smallest atomic unit that needs to be scheduled.
For each task $t \in T$, we define a set of processing modes $M_t \subseteq M$, where exactly one mode must be selected.
A mode represents one possible way to process a task.
Each mode $m \in M$ specifies a processing duration $p_m \geq 0$, a set of required resources $R_m \subseteq R$, and resource demands $q_{mr} \geq 0$ for each resource $r \in R_m \setminus R^{\text{machine}}$.

We classify constraints between tasks into three categories: timing, assignment, and sequencing constraints. 
Each constraint is formally represented as a tuple in a constraint set $C^\text{ConstraintType}$ for a specific constraint type.

\vspace{2 mm}

\noindent \textit{Timing constraints} define temporal relationships between two tasks $i$ and $k$, with an optional delay $l \in \mathbb{Z}$. These are represented by four sets:
\begin{itemize}
\item $\forall(i, k, l) \in C^\text{StartBeforeStart}$: Task $i$ must start before task $k$ starts by at least $l$ time units
\item $\forall(i, k, l) \in C^\text{StartBeforeEnd}$: Task $i$ must start before task $k$ ends by at least $l$ time units
\item $\forall(i, k, l) \in C^\text{EndBeforeStart}$: Task $i$ must end before task $k$ starts by at least $l$ time units
\item $\forall(i, k, l) \in C^\text{EndBeforeEnd}$: Task $i$ must end before task $k$ ends by at least $l$ time units
\end{itemize}

\vspace{2 mm}

\noindent \textit{Assignment constraints} govern resource allocation decisions between tasks. 
Let $m_i$ and $m_k$ be the selected modes for tasks $i$ and $k$, with corresponding resource sets $R_{m_i}$ and $R_{m_k}$. 
Two types of assignment constraints exist:
\begin{itemize}
\item $\forall (i, k) \in C^\text{IdenticalResources}$: Tasks $i$ and $k$ must use the same resources, i.e., $R_{m_i} = R_{m_k}$
\item $\forall (i, k) \in C^\text{DifferentResources}$: Tasks $i$ and $k$ must use different resources, i.e., $R_{m_i} \cap R_{m_k} = \emptyset$
\end{itemize}

\vspace{2 mm}

\noindent \textit{Sequencing constraints} impose restrictions between tasks and apply only when machines are involved.
Let the overlapping machines of tasks $i$ and $k$ be $R^{\text{machine}}_{ik}= \{r \in R^\text{machine} : r \in R_u \land r \in R_v \text{ for some } u\in M_i, v \in M_k\}$.

\begin{itemize}
\item $\forall (i, k) \in C^\text{Consecutive}$: Task $i$ must immediately precede task $k$ for all machines $r \in R^\text{machine}_{ik}$ they are both scheduled on
\item $\forall (i,k,r, l) \in C^\text{SetupTime}$: When task $k$ is scheduled after task $i$ on machine $r \in R^{\text{machine}}_{ik}$, then there is a setup time of $l$ time units
\end{itemize}

\vspace{2 mm}

A feasible solution for the scheduling problem specifies for each task $t \in T$: (i) when it starts and ends and (ii) which processing mode is selected, while respecting all constraints.
The objective is to minimize a weighted sum of common scheduling objective functions.
We define this in detail in the next section.

\subsection{Constraint programming model}
\label{subsec:model}

\subsubsection{Preliminaries.} 
This section provides a brief introduction to important CP concepts related to scheduling. 
It is intended for readers who are familiar with basic optimization concepts, such as variables and constraints, but have limited knowledge of CP. 
A detailed explanation of CP is outside the scope of this paper. 
For a general introduction to CP from an operations research perspective, we refer readers to \citet{baptiste_constraint-based_2001, kanet_constraint_2004,pesant_constraint_2014}.

CP is a paradigm for addressing constraint satisfaction problems. 
A constraint satisfaction problem consists of a finite set of variables, each with a discrete domain, and a set of constraints that must be satisfied. 
CP systematically narrows down the domains of variables through a technique called constraint propagation, which ensures that constraints are effectively communicated across variables.
In addition to constraint propagation, CP also employs search techniques, such as backtracking and large neighborhood search, to systematically explore possible assignments when propagation alone is insufficient to determine a solution.
CP is designed to handle a wide range of (non-linear) constraints, including mathematical, logical, and global constraints.
This flexibility is essential in capturing the unique aspects of scheduling problems, where specialized global constraints such as \NoOverlap~and \Cumulative~are highly effective at domain reduction and compactly formulate the scheduling constraints, while an equivalent MILP formulation would require a large number of linear constraints and big-M reformulations.

A central element of modern CP solvers such as OR-Tools and CP Optimizer is interval variables \citep{laborie_reasoning_2008}.
An interval variable $\nu$ is a special decision variable that is composed of four other decision variables: the start time variable $\nu^\start \geq 0$, the duration variable $\nu^\duration \geq 0$, the end time variable $\nu^\supend \geq 0$, and the presence variable $\nu^\present \in \{0, 1\}$.
Interval variables impose constraints between the start, duration, and end variables depending on their presence status.
If the interval is present ($\nu^\present = 1$), then $\nu^\duration = \nu^\supend - \nu^\start$ is enforced.
When the interval is absent ($\nu^\present = 0$), then there is no such enforcement.
Moreover, specialized scheduling constraints such as \NoOverlap~and \Cumulative~implicitly take into account the presence of interval variables; if an interval is not present, then it is effectively ignored by the constraint.

In CP, there is a close integration between modeling languages and constraint programming solvers.
For instance, CP Optimizer implements the so-called \texttt{Span} constraint which can be used to relate interval variables to each other, whereas OR-Tools does not, resulting in two different solver-specific formulations.
Many recent scheduling studies have used CP Optimizer and present a CP formulation using CP Optimizer-specific syntax from \citet{laborie_ibm_2018}.
As a result, the description of OR-Tools models is underrepresented in the scheduling academic literature, and in the following, we formulate PyJobShop's CP model using OR-Tools syntax. 

\subsubsection{Model formulation.}
This section presents our model's variables, constraints, and objective functions in sequence.

\subsubsection*{Variables.} We introduce the following interval variables:
\begin{itemize}
    \item $\phi_j$: interval variable for each job $j \in J$ with $\phi_j^\present = 1$
    \item $\tau_t$: interval variable for each task $t \in T$ with $\tau_t^\present = 1$
    \item $\mu_m$: interval variable for each mode $m \in M$ with $\mu_m^\duration = p_m$
\end{itemize}
Job ($\phi_j$) and task ($\tau_t$) interval variables are always present, whereas $\mu_m$ can be optional.
The mode duration $\mu_m^\duration$ is set to be fixed to the mode duration, while the task duration and job duration follow from the constraints as defined below.
The relationship between variables is further explained by the constraints in the next paragraphs.

\subsubsection*{Constraints.}
We present the constraints organized by each logical category.
This follows the same structure as in the code implementation.

\vspace{2 mm}

\noindent\textit{Linking jobs to tasks.} 
Job interval variables are not directly scheduled.
Instead, a job starts and ends with its earliest and latest task, respectively.
\begin{subequations}
\label{model:constraints}
\begin{align}
\label{cons:job-span-task-start}
\phi_{j}^\start = \min_{t \in T_{j}} \tau_t^\start &\quad \forall j \in J \\
\label{cons:job-span-task-end}
\phi_{j}^\supend = \max_{t \in T_{j}} \tau_t^\supend &\quad \forall j \in J
\end{align}
Constraints~\eqref{cons:job-span-task-start} set the start time of job $j$ equal to the earliest start time among all its tasks. 
Similarly, Constraints~\eqref{cons:job-span-task-end} ensure that the job's end time corresponds to the completion time of its last task.

\vspace{2 mm}

\noindent\textit{Linking tasks to modes.} 
Task intervals and mode intervals interact in a structured manner. 
Each task requires the selection of exactly one mode. 
The interval variable of a task always starts and ends simultaneously with the interval variables of its modes, regardless of which mode is selected. 
The selected mode determines the duration of the corresponding task interval. 
\begin{align}
\label{cons:select-one-mode}
\sum_{m \in M_t} \mu_m^\present = 1 &\quad \forall t \in T \\
\label{cons:sync-start}
\mu_m^\start = \tau_t^\start &\quad \forall t \in T, m \in M_t \\
\label{cons:sync-end}
\mu_m^\supend = \tau_t^\supend &\quad \forall t \in T, m \in M_t \\
\label{cons:sync-duration}
\mu_m^\present \implies \mu_m^\duration = \tau_t^\duration &\quad \forall t \in T, m \in M_t 
\end{align}
Constraints~\eqref{cons:select-one-mode} ensure that exactly one mode is selected for each task.
Constraints~\eqref{cons:sync-start} and Constraints~\eqref{cons:sync-end} ensure that the mode variable and task interval start and end together, respectively.
Constraints~\eqref{cons:sync-duration} enforce that if the given mode is selected, the duration is synchronized with the corresponding task duration. 
Additionally, the selected mode interval effectively represents the task with a specific resource allocation, as defined next.

\vspace{2 mm}

\noindent\textit{Resource constraints.} 
Each resource type has a specific rule that dictates how many tasks it can process.
Denote $M^R_r = \{m \in M : r \in R_m \}$ as the modes requiring resource $r \in R$.
The following constraints ensure that resource utilization constraints are respected. 

\begin{align}
\label{cons:no-overlap}
\NoOverlap (\{\mu_m : m \in M^R_r \}) &\quad \forall r \in R^{\text{machine}} \\
\label{cons:renewable-resource}
\Cumulative \{\{\mu_m : m \in M^R_r \}, \{q_{mr} : m \in M^R_r \}, Q_r\} &\quad \forall r \in R^{\text{renewable}} \\
\label{cons:non-renewable-resource}
\sum_{m \in M^R_r} \mu_m^\present \cdot q_{mr} \leq Q_r &\quad \forall r \in R^{\text{non-renewable}}
\end{align}
Constraints~\eqref{cons:no-overlap} ensure that mode variables that use a given machine cannot overlap, that is, a machine can only process one task at a time.
Constraints~\eqref{cons:renewable-resource} restrict that mode variables do not exceed the demand of the requested renewable resource at any point in time, while Constraints~\eqref{cons:non-renewable-resource} ensure that the total demand for a non-renewable resource is not exceeded.
As stated before, the global constraints \NoOverlap~and~\Cumulative~explicitly take into account the presence of the interval variables; if a mode interval variable is not present, then it is effectively ignored by the constraint.

\vspace{2 mm}

\noindent\textit{Timing constraints.} Based on the four sets of tuples capturing the timing constraints, the timing constraints are formally defined as follows:
\begin{align}
\label{cons:start-before-start}
\tau_{i}^{\start} + l \le \tau_{k}^{\start} &\quad \forall (i, k, l) \in C^\text{StartBeforeStart} \\
\label{cons:start-before-end}
\tau_{i}^{\start} + l \le \tau_{k}^{\supend} &\quad \forall (i, k, l) \in C^\text{StartBeforeEnd}\\
\label{cons:end-before-start}
\tau_{i}^{\supend} + l \le \tau_{k}^{\start} &\quad \forall (i, k, l) \in C^\text{EndBeforeStart}\\
\label{cons:end-before-end}
\tau_{i}^{\supend} + l \le \tau_{k}^{\supend} &\quad \forall (i, k, l) \in C^\text{EndBeforeEnd}
\end{align}
Each of the constraints in \eqref{cons:start-before-start}--\eqref{cons:end-before-end} restricts the timing between pairs of task variables.

\vspace{2 mm}

\noindent\textit{Assignment constraints.} 
Some pairs of tasks must use either identical or different resources. 
For any such pair of tasks $i$ and $k$, we need to ensure their selected modes are compatible. 
\begin{align}
    \label{cons:identical-resources}
    \mu_{m_i}^\present \leq \sum_{\substack{m_k \in M_k \\ \text{s.t. } R_{m_k} = R_{m_i}}} \mu_{m_k}^\present &\quad \forall (i, k) \in C^{\text{IdenticalResources}}, m_i \in M_i \\ 
    \label{cons:different-resources}
    \mu_{m_i}^\present \leq \sum_{\substack{m_k \in M_k \\ \text{s.t. } R_{m_k} \cap R_{m_i} = \emptyset}} \mu_{m_k}^\present &\quad \forall (i, k) \in C^{\text{DifferentResources}}, m_i \in M_i
\end{align}
Constraints~\eqref{cons:identical-resources} ensure that if mode $m_i$ is selected for task $i$, then at least one mode with identical resources must be selected for task $k$. 
Similarly, Constraints~\eqref{cons:different-resources} ensure that if mode $m_i$ is selected for task $i$, then at least one mode with disjoint resources must be selected for task $k$.

\vspace{2 mm}

\noindent\textit{Sequencing constraints.}
Setting up sequencing constraints in OR-Tools requires more setup because OR-Tools does not provide an interface for sequencing variables like CP Optimizer \citep{laborie_ibm_2018}.
Instead, with OR-Tools, a sequence of intervals can be represented by a complete graph combined with the global \texttt{Circuit} constraint to select a specific ordering of the intervals.
For each machine $r \in R^\text{machine}$, define a complete graph where the node set $V_r = M^R_r \cup \{0\}$ includes all modes that require machine $r$ plus a dummy node $0$. 
The arcs consist of all possible node pairs, including self-loops. 
We introduce binary variables $B_r = \{b_{ruv} \in \{0, 1\} : u,v \in V_r\}$, where $b_{ruv} = 1$ indicates that arc $(u,v)$ is selected in the graph of machine $r$, and $0$ otherwise. 
Let $t_{m}$ denote the task associated with mode $m$.
Then, the sequencing constraints are specified as follows.
\begin{align}
    \label{cons:circuit}
    &\texttt{Circuit($B_r$)} &\quad \forall r \in R^{\text{machine}} \\
    \label{cons:arcs-implies-presence}
    &b_{ruv} \implies \mu_{u}^\present \land \mu_{v}^\present &\quad \forall u, v \in M^R_r, r \in R^{\text{machine}} \\
    \label{cons:self-arc}
    &b_{ruu} \implies \neg \mu_{u}^\present &\quad \forall u \in M^R_r, r \in R^{\text{machine}} \\ 
    \label{cons:dummy-self-arc}
    &b_{r00} \implies \neg \mu_{u}^\present &\quad \forall u \in M^R_r, r \in R^{\text{machine}} \\ 
    \label{cons:setup-time}
    &b_{ruv} \implies \mu_{u}^\supend + s_{t_u, t_v, r} \leq \mu_{v}^\start &\quad \forall u, v \in M^R_r, r \in R^{\text{machine}}
\end{align}
Constraints~\eqref{cons:circuit} ensure that the selected arcs form a single sub-tour, starting and ending at the dummy node, thus establishing an order for the mode intervals on machine $r$.
Constraints~\eqref{cons:arcs-implies-presence} guarantee that if an arc is selected, the corresponding mode intervals must be present. 
For modes not part of the selected sub-tour, self-loops must be selected, and Constraints~\eqref{cons:self-arc} ensure that the corresponding mode interval is absent.
Constraints~\eqref{cons:dummy-self-arc} address the special case where if the dummy self-arc is selected, all mode intervals must be absent. 
This occurs when a resource is not allocated any tasks.
Constraints~\eqref{cons:setup-time} then ensure that there is an end-before-start relationship including setup times. 
If $(t_u,t_v,r, l) \in C^{\text{SetupTime}}$, then the setup time $s_{t_u, t_v, r}$ is given by $l$, otherwise it is zero.

The consecutive constraints are defined as follows. 
\begin{align}
    \mu_{u}^\present \land \mu_{v}^\present \implies b_{ruv} &\quad \forall (i, k) \in C^{\text{Consecutive}}, r \in R^{\text{machine}}_{ik}, u \in M_i \cap M^R_r, v \in M_k \cap M^R_r
\end{align}
This constraint ensures that if both modes $u$ and $v$, belonging to tasks $i$ and $k$ respectively, are present, then it must be part of the sub-tour selected by the \texttt{Circuit} constraint.
Consequently, this ensures that tasks $i$ and $k$ are scheduled consecutively on the set of machines that they are both allocated to.

\end{subequations}

\subsubsection*{Objective function.}
PyJobShop supports the following common objective functions:
\begin{itemize}
    \item \text{Makespan}: $\max_{t \in T} \tau_t^\supend$
    \item \text{Total weighted flow time}: $\sum_{j \in J} w_j (\phi_j^\supend - r_j)$
    \item \text{Total weighted tardiness}: $\sum_{j \in J} w_j \max ( \phi_j^\supend - d_j, 0 ) $
    \item \text{Total weighted earliness}: $\sum_{j \in J} w_j \max ( d_j - \phi_j^\supend, 0 )$
    \item \text{Total weighted number of tardy jobs}: $\sum_{j \in J} w_j \mathbbm{1}\{ \phi_j^\supend > d_j\}$
    \item \text{Maximum tardiness}: $\max_{j \in J} w_j \max{(\phi_j^\supend - d_j, 0)}$
    \item \text{Maximum lateness}: $\max_{j \in J} w_j (\phi_j^\supend - d_j)$ 
\end{itemize}
Let $F$ denote the set of objective functions and let $w^f$ denote the weight of objective $f \in F$.
The overall objective is to find a solution $x$ in the set of feasible solutions $X$ that minimizes the weighted sum of all objective functions:
\begin{align}
    \min_{x \in X} \sum_{f \in F} w^f \cdot f(x)
\end{align}

%% file: sections/variants.tex
\section{Supported scheduling problems}
\label{sec:variants}

PyJobShop's scheduling model, introduced in Section~\ref{sec:problem}, supports a wide range of scheduling problems. 
In this section, we describe several common scheduling variants that can be modeled and solved with PyJobShop. 
We discuss machine scheduling variants in Section~\ref{subsec:machine} and project scheduling variants in Section~\ref{subsec:project}.

\subsection{Machine scheduling} \label{subsec:machine}
To demonstrate the supported machine scheduling variants, we rely on Graham's notation $\mbox{$\alpha \mid \beta \mid \gamma$ }$, which is widely used in the scheduling community to classify scheduling problem models \citep{graham_optimization_1979}.
The first field ($\alpha$) specifies the machine environment, the second field ($\beta$) specifies the constraints, and the third field ($\gamma$) specifies the objective.
Table~\ref{table:classification} presents the most common scheduling characteristics for each field based on the classification schemes described from~\citet{framinan_manufacturing_2014,pinedo_scheduling_2016,framinan_deterministic_2019,dauzere-peres_flexible_2024}.

In terms of machine environments, all classic environments (1, P, F, HF, O, J, FJ; see Table~\ref{table:classification} for their meaning) are supported since most of them can be cast as a variant of the flexible job shop environment.
Beyond the classic environments, an environment commonly found in realistic manufacturing settings is the assembly scheduling problem ($\circ \rightarrow \circ$), which involves the scheduling of tasks in a concurrent fashion \citep{framinan_deterministic_2019}.
Assembly scheduling generalizes the order scheduling environment \citep{leung_order_2005}, where a job can consist of multiple tasks with their own processing sequence that need to be completed before it is considered completed.
For example, a common manufacturing environment is $Pm \rightarrow 1$, which consists of parallel machines in the first stage that produce parts of some product, followed by a single assembly line that can only start when all individual parts have been produced.

\citet{dauzere-peres_flexible_2024} extends the environment field with multi-mode (MM), flexible sequencing (FS), multi-resource (MR), flexible processing planning (FP), and distributed (D).
Multiple modes (MM) are supported and described as a core part of our model.
Flexible sequencing (FS) allows interchanging the order in which a job's tasks are performed, for instance, processing tasks 1, 2, 3 or 1, 3, 2.
Additionally, some tasks may not be processed simultaneously (e.g., tasks 2 and 3), which can achieved by introducing a fictitious machine and modes to prevent overlap \citep{dauzere-peres_flexible_2024}.
Multi-resource (MR) is a generalization of multiple modes, where each mode describes a set of required skills that must be satisfied by any resource that possesses that skill.
A more fitting name for multi-resource is multi-skilled, which is the terminology used in the project scheduling literature \citep{snauwaert_classification_2023}. 
Multi-resource is not supported in PyJobShop but can be implemented as a multi-mode problem if the number of resource-skill combinations is limited.
Flexible processing planning (FP) and distributed (D) scheduling define alternative routes for a job's tasks, meaning that there are multiple distinct ways to process a job, e.g., either by processing tasks 1 and 2 or by processing tasks 3 and 4.
Both these extensions are not supported because they require tasks to be optional.

In terms of constraints, PyJobShop supports release dates ($r_j$) and deadlines ($\overline{d}_j$), which restrict job start and completion times, as well as due dates ($d_j$) for tardiness-based objectives.
Capacity-based resources (res) describe the use of renewable resources in machine scheduling problems (as opposed to only machines) and are explicitly part of our scheduling model.
Bills of materials (bom), also known as arbitrary precedence graphs \citep{kasapidis_flexible_2021}, are supported since timing constraints can be arbitrarily imposed between tasks.
Sequence-dependent setup times ($s_{ijk}$) are supported as sequencing constraints.
Blocking (block) requires that a task can occupy a resource longer than its processing duration requires because its successor tasks can not yet start.
This is supported by allowing variable task durations and generalized precedence constraints (end-at-start).
Buffers ($b_{i}$) can be modeled by introducing additional resources for each machine, each with a capacity and precedence constraints with blocking.
No waiting time (no-wait) between consecutive tasks is possible by combining end-before-start and start-before-end constraints.
Task overlap (overlap) between tasks is possible as nothing restricts tasks from overlapping in a resource setting, only the resource capacity (cumulative) or one at a time (machine).
Basic forms of machine breakdowns or unavailability (brkdwn) can be added by introducing dummy tasks that occupy a machine for a fixed period of time.
General time lags ($d_{ii'}^{kk'}$) incorporate extra delays to the precedence constraints, which is supported if the time lag is machine-independent.

We have intentionally chosen not to support permutation constraints (prmu) due to the complexity of implementing them robustly alongside other features. 
Permutation constraints require a specific mapping between tasks and machines to ensure correctness, which is not feasible with the current interface. 
Additionally, we do not support no-idle constraints (no-idle), which require machines to continuously process tasks without any idle time between them. 
Batching problems (p-batch) are also unsupported, as they involve an additional decision regarding which tasks should be processed simultaneously. 
Similarly, pre-emption (prmp) is not supported, as it requires making another decision about how to split tasks.

PyJobShop supports the most common objective functions, including makespan ($C_{\max}$), total weighted flow time (TWFT), total weighted tardiness (TWT), total weighted earliness (TWE), the weighted number of tardy jobs (TWNTJ), maximum lateness ($L_{\max}$), and maximum tardiness ($T_{\max}$).
These objective functions share a common feature: they are all based on the completion times of jobs or tasks.
There is currently no support for objectives based on other factors, such as resource workload.
For a detailed overview of relevant objectives in machine scheduling settings, we refer to \citet{ostermeier_review_2023}.

\begin{table}
\caption{%
  Classification scheme for machine scheduling problems. Supported features by PyJobShop are indicated with a black square, while unsupported features are indicated with a white square.
}
\label{table:classification}
{ 
\begin{multicols}{3}
\footnotesize
\setlength{\parindent}{0pt}
\setlist[itemize]{left=0pt}

\textbf{Machine environments (\bm{$\alpha$})}
\begin{itemize}[label=$\square$,itemsep=0pt]
  \item[$\blacksquare$] $1$: single machine
  \item[$\blacksquare$] P: parallel machines
  \item[$\blacksquare$] F: flow shop
  \item[$\blacksquare$] HF: hybrid flow shop
  \item[$\blacksquare$] O: open shop
  \item[$\blacksquare$] J: job shop
  \item[$\blacksquare$] FJ: flexible job shop
  \item[$\blacksquare$] $\circ \rightarrow \circ$: assembly scheduling
  \item[$\blacksquare$] MM: multi-mode
  \item[$\blacksquare$] FS: flexible sequencing
  \item MR: multi-resource
  \item FP: flexible processing planning
  \item D: distributed
\end{itemize}

\columnbreak
\textbf{Constraints (\bm{$\beta$})}
\begin{itemize}[label=$\square$,itemsep=0pt]
    \item[$\blacksquare$] $r_j$: release dates
    \item[$\blacksquare$] $d_j$: due dates
    \item[$\blacksquare$] $\overline{d}_j$: deadlines
    \item[$\blacksquare$] res: capacity-based resources
    \item[$\blacksquare$] $s_{ijk}$: setup times
    \item[$\blacksquare$] bom: bills of materials
    \item[$\blacksquare$] \text{block}: jobs are blocked
    \item[$\blacksquare$] \text{$b_i$}: buffers
    \item[$\blacksquare$] \text{no-wait}: jobs may not wait
    \item[$\blacksquare$] overlap: overlapping tasks
    \item[$\blacksquare$] \text{brkdwn}: breakdowns
    \item[$\blacksquare$] $d_{ii'}^{kk'}$: general time lags
    \item no-idle: machine cannot be idle
    \item \text{prmu}: permutation constraint
    \item $p$-batch: simultaneous processing
    \item \text{prmp}: pre-emption


\end{itemize}

\columnbreak
\textbf{Objectives ($\bm{\gamma}$)}

\begin{itemize}[label=$\square$,itemsep=0pt]
    \item[$\blacksquare$] $C_{\max}$: makespan
    \item[$\blacksquare$] TWFT: total weighted flow time
    \item[$\blacksquare$] TWT: total weighted tardiness
    \item[$\blacksquare$] TWE: total weighted earliness
    \item[$\blacksquare$] TWNTJ: total weighted number of tardy jobs
    \item[$\blacksquare$] $T_{\max}$: maximum tardiness
    \item[$\blacksquare$] $L_{\max}$: maximum lateness
    \item $W_T$: total workload of machines
\end{itemize}

\end{multicols}
}
\end{table}

\subsection{Project scheduling} \label{subsec:project}
The original scope of PyJobShop was to specifically solve machine scheduling problems. 
However, because we introduced concepts that are also prevalent in the project scheduling literature, such as modes and capacity-based resources, there are many project scheduling variants that PyJobShop can also solve.

Project scheduling and machine scheduling share many similar ideas, but the nomenclature is slightly different \citep{demeulemeester_project_2006}.
Jobs are called projects, although it is common to only have one project, and tasks are referred to as activities and events.
Resources include renewable and non-renewable resources, or double-constrained resources, but also cumulative and partially renewable resources, the latter of the two which are not supported in PyJobShop.
Machines are uncommon in the project scheduling literature, but there are variants, such as the multi-skilled project scheduling problem, that use disjunctive resources \citep{snauwaert_classification_2023}.

We outline here the most well-known project scheduling variants that are supported by PyJobShop.
We refer to the surveys of \cite{hartmann_updated_2022} and \cite{gomez_sanchez_resource-constrained_2023} for an extensive survey of the project scheduling variants.
The project scheduling problem is the most basic variant, which considers tasks without any resource constraints, but tasks are constrained by precedence relationships.
The resource-constrained project scheduling problem (RCPSP) introduces finite resources, each task requiring a subset of these resources for a fixed duration and resource requirement.
The multi-mode variant of RCPSP (MMRCPSP) extends the RCPSP by allowing tasks to be executed in various modes, each with distinct duration and resource requirements, adding flexibility that is similar to the extension from job shops to flexible job shops.
Another commonly studied variant of the RCPSP is the one with generalized precedence constraints (RCPSP/max) \citep{schutt_scheduling_2013}. 
Generalized precedence relations express relations of start-to-start, start-to-end, end-to-start, and end-to-end between pairs of tasks, as described earlier in our scheduling model.
Finally, the resource-constrained multi-project scheduling problem (RCMPSP) introduces multiple concurrent projects, each with a set of tasks that need to be completed.
The goal of all these problem variants is to minimize the makespan. 

%% file: sections/software.tex
\section{Software}
\label{sec:software}
The PyJobShop package is developed in a GitHub repository located at \url{https://github.com/PyJobShop/PyJobShop}. 
This repository contains the source code, including tests, documentation, and examples that introduce new users to PyJobShop.
The documentation of PyJobShop is hosted at \url{https://pyjobshop.org}.
Users can directly install PyJobShop from the Python package index using \texttt{pip install pyjobshop}, which comes with OR-Tools by default.
Additional installation instructions for CP Optimizer are provided in the documentation.

PyJobShop borrows many ideas from PyVRP~\citep{wouda_pyvrp_2024}, an open-source solver for vehicle routing problems.
In particular, the simple modeling framework, extensive documentation, and numerous examples have been shown to bring much value to users of PyVRP from academia as well as industry.

In the following sections, we elaborate on the package structure, provide a modeling example and describe how to extend the package, respectively.

\subsection{Package structure}
The top-level package of the \texttt{pyjobshop} namespace contains most user components.
\begin{itemize}
    \item \texttt{ProblemData.py}: Contains the \texttt{ProblemData} class, defining the problem data instance to be solved, as well as the \texttt{Job}, \texttt{Task}, \texttt{Machine}, \texttt{Renewable}, \texttt{NonRenewable}, \texttt{Mode} and \texttt{Constraints} classes.
    \item \texttt{Model.py}: The modeling interface to build a \texttt{ProblemData} instance step-by-step.
    \item \texttt{Solution.py}: The class that describes a solution.
    \item \texttt{Result.py}: The result of a solver run, including the best-found solution along with solver statistics.
    \item \texttt{read.py}: A function to read a variety of scheduling benchmark instances.
    \item \texttt{solve.py}: A dedicated solving function.
    \item \texttt{cli.py}: The command-line interface, mostly for internal use and benchmarking.
    \item \texttt{solvers/}: The solver module implements the scheduling model using constraint programming solvers, \texttt{ortools} and \texttt{cpoptimizer}. Each solver is encapsulated in a solver-specific \texttt{Solver} class, which manages implementations of \texttt{Variables}, \texttt{Constraints}, and \texttt{Objective} classes.
\end{itemize}

\subsection{Example use}
\begin{figure}
    \caption{Gantt chart produced by the code from Listing 1. Each bar represents one task and the colors depict the job it belongs to.}
    \label{fig:listing1}
    \centering
    \includegraphics[width=\linewidth]{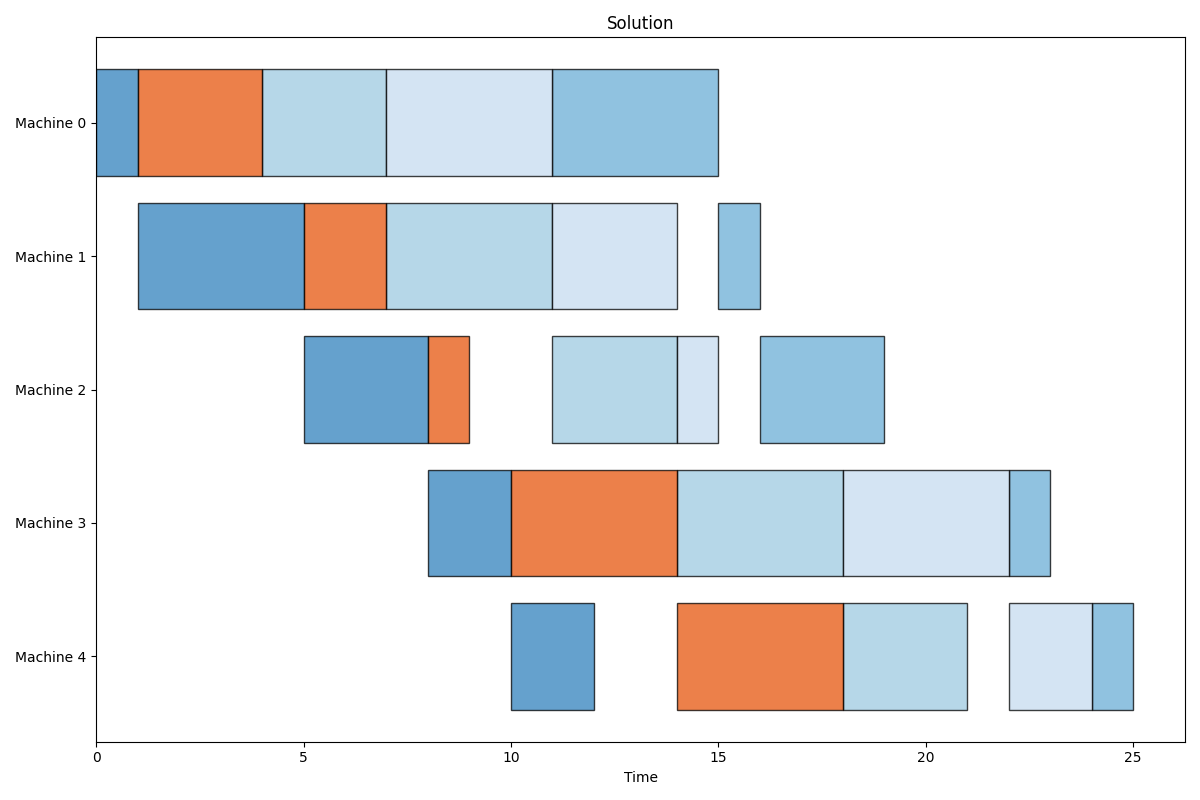}
\end{figure}

The primary interface for PyJobShop is the \texttt{Model} class, which provides an intuitive domain-specific interface for scheduling problems.
Listing 1 presents a complete example of modeling and solving a flow shop problem.
Users first create instance components (machines, jobs, tasks, modes, constraints) through dedicated model methods and then invoke the \texttt{solve} method to obtain a result object containing the solution and solver statistics.
The solution can then be plotted using \texttt{plot\_machine\_gantt} from the \texttt{pyjobshop.plot} module, which plots a Gantt chart as shown in Figure~\ref{fig:listing1}.
We refer to the documentation for more examples, which include more examples from machine scheduling as well as project scheduling.

\begin{lstlisting}[
float,
style=mypython, 
language=Python, 
caption=Modeling a flow shop problem.,
label=lst:machine]
import random
from pyjobshop import Model
from pyjobshop.plot import plot_machine_gantt

random.seed(42)

model = Model()
machines = [model.add_machine() for idx in range(5)]

for job_idx in range(5):
    job = model.add_job()
    tasks = [model.add_task(job=job) for idx in range(5)]

    for idx in range(len(tasks)):
        task = tasks[idx]
        machine = machines[idx]
        processing_time = random.randint(1, 5)
        model.add_mode(task, machine, processing_time)

    for idx in range(len(tasks) - 1):
        pred = tasks[idx]
        succ = tasks[idx + 1]
        model.add_end_before_start(pred, succ)

result = model.solve(time_limit=10)
data = model.data()
solution = result.best

plot_machine_gantt(solution, data)
\end{lstlisting}

\subsection{Extending PyJobShop}
While PyJobShop implements common constraints for scheduling problems, some users may require more specialized constraints.
Users can easily extend the PyJobShop framework to suit their needs.
Extensions often fall into one of the following categories: adding new constraints, adding new objective functions, or adding new decision variables.
A common workflow for making changes is by modifying, in order, the following classes: \texttt{ProblemData} (how to represent the new feature as data), \texttt{Model} (how to interact with the new feature through the user interface) and the corresponding CP solver's \texttt{Solver} class (how to implement the new feature in terms of CP variables and constraints).

We are open to new contributions and encourage users to contribute with new features, examples, and documentation improvements.
To discuss potential extensions or improvements, we recommend first opening an issue on the GitHub repository.

%% file: sections/experiments.tex
\section{Numerical experiments}
\label{sec:experiments}
PyJobShop allows for easy modeling of many scheduling problems through one single interface. 
We leverage this to evaluate and compare the performance of OR-Tools and CP Optimizer on a large variety of scheduling instances.
Section~\ref{subsec:instances} covers the benchmark instances and Section~\ref{subsec:computational} describes the computational aspects. 
The results are presented in Section~\ref{subsec:results}, and we discuss the results in Section~\ref{subsec:discussion}.

\subsection{Benchmark instances}
\label{subsec:instances}
We conducted our experiments on 9,280 instances total, both from machine scheduling and project scheduling instances.
Table~\ref{table:machine-instance-stats} summarizes the instance characteristics.
In the following, we describe in more detail the choices for these instances, discussing machine and project scheduling separately.

{
\begin{table}
\caption{Summary of instance statistics of benchmark instances used for the numerical experiments.
The table reports the minimum, average, and, maximum number of tasks and resources over all instances for the given statistic.}
\label{table:machine-instance-stats}
\renewcommand{\arraystretch}{1.25}
\centering
\input{tables/machine-instance-stats}
\end{table}
}

\paragraph{Machine scheduling instances.}
The first set covers a subset of machine scheduling instances used in \citet{naderi_mixed-integer_2023}.
We consider nine different problem variants, including the job shop problem (JSP), flexible job shop problem (FJSP), no-wait permutation flow shop problem (NW-PFSP), non-permutation flow shop problem (NPFSP), hybrid flow shop problem (HFSP), permutation flow shop problem (PFSP), sequence-dependent setup times PFSP (SDST-PFSP), total completion time PFSP (TCT-PFSP) and total tardiness PFSP (TT-PFSP).
We refer to \citet{naderi_mixed-integer_2023} for a precise definition of these problem variants and benchmark instances.

Four out of nine selected variants require explicit permutation constraints (PFSP, SDST-PFSP, TCT-PFSP, TT-PFSP), which are not directly supported by PyJobShop.
Nevertheless, to evaluate the performance between the two CP solvers and to replicate the study performed by \citet{naderi_mixed-integer_2023}, we have implemented permutation constraints specifically for this experiment.
However, OR-Tools lacks an efficient way of implementing the permutation constraint, requiring the basic sequencing constraints from Section~\ref{subsec:model} that result in a large number of variables and constraints. 
This limitation forced us to restrict the instances to at most 100 jobs, as larger problems were computationally intractable. 
Finally, we note that the NW-PFSP does not require explicit permutation constraints as this can be implemented using end-to-start constraints, which implies a permutation in the context of flow shops.

We excluded three problem variants from \citet{naderi_mixed-integer_2023}: (i) the open shop problem because those instances are all trivially solved, (ii) the parallel machines problem because this problem admits much more efficient scheduling models than the one used in PyJobShop, and (iii) the distributed PFSP, which is not supported even when implementing the permutation constraint.

For all machine scheduling instances, we use the best-known solutions from the data provided in \cite{naderi_mixed-integer_2023}.
The best-known solutions may already be outdated since those experiments were conducted in 2021, though it is outside the scope of this paper to find all the most recent best-known solutions.

\paragraph{Project scheduling instances.}
From the project scheduling literature, we included instances from the resource-constrained project scheduling problem (RCPSP), the multi-mode project scheduling problem (MMRCPSP), and the resource-constrained multi-project scheduling problem (RCMPSP).
These instances differ from machine scheduling instances because they require renewable and non-renewable resources instead of machines.

For RCPSP, we use instances from two well-established benchmarks: PSPLIB \citep{kolisch_psplib_1997} with 30, 60, 90, and 120 tasks, and RG300 \citep{debels_decomposition-based_2007} with 300 tasks.
For MMRCPSP, we use instances from MMLIB~\citep{van_peteghem_experimental_2014}, specifically, the instance sets MMLIB50 and MMLIB100 containing instances with 50 and 100 tasks, respectively.
For RCMPSP, we use instances from MPLIB1~\citep{van_eynde_resource-constrained_2020}, specifically, instance set 3, which has the largest number of tasks (1488) of the MPLIB1 dataset.
The most recent best-known solutions are taken from \citet{operations_research__scheduling_research_group_rcpsp_2025} for RCPSP and MMRCPSP, and taken from \citet{bredael_multi-project_2023} for RCMPSP.

\subsection{Computational details}\label{subsec:computational}
Each instance was solved using eight cores of an AMD EPYC 9654 CPU with a 900-second time limit. 
For each instance, we compute the optimality gap and the relative percentage deviation (RPD):
\begin{align}
  \text{Gap} = \frac{\text{UB} - \text{LB}}{\text{UB}} \times 100  \quad\quad& \text{RPD} = \frac{\text{UB} - \text{BKS}}{\text{BKS}} \times 100
\end{align}
where UB and LB are the obtained upper and lower bounds, respectively, and BKS is the best-known solution value.
Note that for some TT-PFSP instances, the UB or BKS is 0. 
In such cases, if the numerator is nonzero, we set the corresponding metric to 100, otherwise, we set it to 0.

We used OR-Tools v9.11.4210 and CP Optimizer v22.1.1.0, each with their respective implementations of the base scheduling model presented in Section~\ref{sec:problem} but slightly modified to accommodate permutation constraints. 
All code, data, and results are available in our GitHub repository at \url{https://github.com/PyJobShop/Experiments}.

\subsection{Results}
\label{subsec:results}

{
\begin{table}
\caption{%
RPD and optimality gap averaged over all feasibly solved instances per problem variant, comparing OR-Tools and CP Optimizer with a 900-second time limit.
}
\label{table:gap-per-problem}
\renewcommand{\arraystretch}{1.4}
\centering
\input{tables/gap-per-problem}

\end{table}
}

Table~\ref{table:gap-per-problem} reports the average RPD and optimality gap across all problem categories by solver.
Instances with infeasible solutions were excluded; however, the majority were solved feasibly, with the exception being CP Optimizer unable to find feasible solutions to 4\% of the MMRCPSP instances.
Overall, the results show that OR-Tools is highly competitive with CP Optimizer. 
It obtains comparable solutions to the JSP, FJSP, NW-PFSP, and all project scheduling variants, even obtaining slightly lower average RPDs on the FJSP, NW-PFSP, and MM-RCPSP. 
Additionally, OR-Tools consistently obtains better optimality gaps than CP Optimizer on all but two permutation problems, demonstrating OR-Tool's ability to compute stronger lower bounds than CP Optimizer.

CP Optimizer clearly excels in solving permutation scheduling problems, even after limiting instance sizes to accommodate OR-Tools. 
Expressing sequencing constraints in OR-Tools is less efficient than in CP Optimizer and this significantly impacts OR-Tools' performance on these problems. 
Moreover, CP Optimizer generally handles large-scale instances more efficiently. 
For instance, CP Optimizer outperforms OR-Tools on NPFSP and HFSP, which include instances with up to 48,000 and 2,000 tasks (and up to 10,000 modes), respectively. 
A notable exception is the NW-PFSP with up to 48,000 tasks, where OR-Tools unexpectedly outperforms CP Optimizer.
After our benchmark, we found that CP Optimizer's performance on NW-PFSP strongly benefits from adding permutation constraints, and based on the results in \citet{naderi_mixed-integer_2023}, it should be expected that CP Optimizer achieves an average RPD of roughly 1.6\% with a 3,600-second time limit, outperforming OR-Tools (3.47\%).

For project scheduling problems, OR-Tools and CP Optimizer perform similarly, with RPD differences of less than 0.5\% across all variants. 
CP Optimizer performs better on RCPSP and RCMPSP, whereas OR-Tools achieves better results on MMRCPSP. 
Notably, both solvers combined found over 22, 13, and 2,025 new best-known solutions to RCPSP, MMRCPSP, and RCMPSP instances, respectively, further strengthening the case of using CP for solving project scheduling problems.

\subsection{Discussion}
\label{subsec:discussion}
Our results establish OR-Tools as a strong alternative to CP Optimizer for solving various scheduling problems.
This contradicts the conclusion of \citet{naderi_mixed-integer_2023}, who claim in their supplementary material (p. 9) that “OR-Tools does not yield a performance that is remotely close to that of CP Optimizer.” 
In the main paper (footnote 7), they further suggest that OR-Tools struggles with assignment variables compared to CP Optimizer based on their results for the FJSP and HFSP.
We reviewed the OR-Tools implementation used in their experiments and found that their models for FJSP and HFSP were inefficiently formulated (see Appendix~\ref{app:fjsp_naderi} for details).
When we used their implementation on the same set of FJSP instances and identical computational setup as our numerical experiments, we obtained an average RPD of 7.55\%.
In contrast, our optimized OR-Tools model achieved an average RPD of just 0.68\%.
This suggests that the main reason for their observed inferior performance of OR-Tools was primarily due to suboptimal modeling choices, though the use of an older version of the CP-SAT solver may have also played a role. 

Despite OR-Tools’ strong results, CP Optimizer still holds advantages in specific areas. 
It scales better to larger instances and performs well on permutation-based scheduling problems, and from our personal experience, CP Optimizer often finds high-quality solutions more quickly. 
The advantage of CP Optimizer in large-scale instances could be attributed to its iterative diving search method, which aggressively dives into the search tree without backtracking \citep{laborie_ibm_2018}.
This is also quantified by \citet{da_col_industrial-size_2022}, where the authors demonstrate that CP Optimizer outperformed OR-Tools by 6--42\% on job shop scheduling problems with 10,000--100,000 tasks, and successfully solved instances with up to 1,000,000 tasks where OR-Tools failed to find solutions.
We remark that \citet{da_col_industrial-size_2022} only studied the performance of both solvers with at most four cores, whereas the recommended number of cores for OR-Tools is at least eight \citep{perron_cp-sat_2024}.

When looking at the broader picture, both solvers demonstrate the effectiveness of CP for scheduling problems. 
As shown by \citet{naderi_mixed-integer_2023}, CP Optimizer generally outperforms MILP solvers in machine scheduling, and our results for CP Optimizer are largely consistent with theirs.
While CP sometimes yields relatively high RPDs, such as around 10\% for NPFSP, this is expected since it is compared to solutions from specialized (meta)heuristics which are fine-tuned for these problems.
However, the key advantage of CP lies in its flexibility as it can easily accommodate new constraints that arise in real-world applications, whereas specialized heuristics often require extensive modifications for each new feature.

%% file: tables/machine-instance-stats.tex
\renewcommand{\arraystretch}{1.4}
\begin{tabular}{c@{\hskip 25pt}lc@{\hskip 10pt}rrr@{\hskip 10pt}rrr}
\toprule
& & & \multicolumn{3}{c}{\# Tasks} & \multicolumn{3}{c}{\# Resources} \\
\cmidrule(lr){4-6} \cmidrule(lr){7-9}
& Problem & \# Instances & Min. & Avg. & Max. & Min. & Avg. & Max. \\
\midrule
\multirow{5}{*}{\rotatebox[origin=c]{90}{Non-permutation}} 
& JSP & 242 & 36 & 511 & 2000 & 5 & 15 & 20 \\
& FJSP & 289 & 12 & 322 & 1477 & 4 & 11 & 20 \\
& NW-PFSP & 360 & 100 & 12610 & 48000 & 5 & 31 & 60 \\
& NPFSP & 360 & 100 & 12610 & 48000 & 5 & 31 & 60 \\
& HFSP & 1440 & 250 & 938 & 2000 & 15 & 30 & 50 \\
\midrule
\multirow{4}{*}{\rotatebox[origin=c]{90}{Permutation}}
& PFSP & 120 & 100 & 1496 & 6000 & 5 & 19 & 60 \\
& SDST-PFSP & 360 & 100 & 661 & 2000 & 5 & 12 & 20 \\
& TCT-PFSP & 120 & 100 & 1496 & 6000 & 5 & 19 & 60 \\
& TT-PFSP & 135 & 500 & 1500 & 2500 & 10 & 30 & 50 \\
\midrule
\multirow{3}{*}{\rotatebox[origin=c]{90}{Project}}
& RCPSP & 2520 & 32 & 122 & 302 & 3 & 4 & 4 \\
& MMRCPSP & 1080 & 52 & 77 & 102 & 4 & 4 & 4 \\
& RCMPSP & 2254 & 1488 & 1488 & 1488 & 4 & 4 & 4 \\
\bottomrule
\end{tabular}

%% file: tables/gap-per-problem.tex
\begin{tabular}{l@{\hskip 25pt}lcc@{\hskip 10pt}cc}
\toprule
& & \multicolumn{2}{c}{RPD (\%)} & \multicolumn{2}{c}{Gap (\%)}  \\
\cmidrule[0.4pt](lr){3-4} \cmidrule[0.4pt](lr){5-6}
& Problem & {OR-Tools} & {CP Optimizer} & {OR-Tools} & {CP Optimizer} \\
\midrule
\multirow{6}{*}{\rotatebox[origin=c]{90}{Non-permutation}} 
& JSP & 1.98 & \textbf{1.80} & \textbf{3.40} & 4.09 \\
& FJSP & \textbf{0.68} & 0.99 & \textbf{1.04} & 27.67 \\
& NW-PFSP & \textbf{3.47} & 7.18 & \textbf{50.51} & 57.87 \\
& NPFSP & 13.53 & \textbf{8.88} & \textbf{16.29} & 25.58 \\
& HFSP & 13.34 & \textbf{6.58} & \textbf{11.94} & 66.54 \\
\cmidrule(r){2-6} 
& \textit{Average} & 6.60  &  \textbf{5.08}  &  \textbf{16.63}  &  36.35 \\
\midrule
\multirow{5}{*}{\rotatebox[origin=c]{90}{Permutation}} 
& PFSP & 7.49 & \textbf{2.54} & 10.61 & \textbf{7.03} \\
& SDST-PFSP & 8.82 & \textbf{4.41} & 30.75 & \textbf{28.24} \\
& TCT-PFSP & 10.31 & \textbf{3.31} & \textbf{21.48} & 26.69 \\
& TT-PFSP & 53.14 & \textbf{20.79} & \textbf{66.86} & 72.43 \\

\cmidrule(r){2-6} 
& \textit{Average} & 19.94  &  \textbf{7.76}  &  \textbf{32.43}  &  33.6 \\
\midrule
\multirow{4}{*}{\rotatebox[origin=c]{90}{Project}}
& RCPSP & 0.93 & \textbf{0.46} & \textbf{3.46} & 4.26 \\
& MMRCPSP & \textbf{0.18} & 0.27 & \textbf{0.94} & 6.65 \\
& RCMPSP & -0.52 & \textbf{-0.85} & \textbf{14.91} & 67.36 \\
\cmidrule(r){2-6} 
& \textit{Average} & 0.20  &  \textbf{-0.04}  &  \textbf{6.44}  &  26.09 \\
\bottomrule
\end{tabular}

%% file: sections/conclusion.tex
\section{Conclusion}
\label{sec:conclusion}
This paper introduced PyJobShop, an open-source Python library for solving scheduling problems with constraint programming.
PyJobShop provides an easy-to-use modeling interface that allows users to solve scheduling problems without having to know the details of constraint programming.
We used PyJobShop to conduct large-scale numerical experiments on a wide variety of scheduling problems from the machine scheduling and project scheduling literature.
Our findings show that OR-Tools is highly competitive with CP Optimizer, particularly in job shop and project scheduling problems, where it often matches and sometimes even surpasses CP Optimizer’s performance.
CP Optimizer excels in permutation scheduling problems, benefiting from its ability to handle sequencing constraints more efficiently, and large-scale scheduling problems.

For future work, we see three interesting directions.
First, there are many more scheduling variants that can be supported by PyJobShop.
This includes task selection problems, which include the decision whether to schedule a task or not \citep{kis_job-shop_2003} and would support the modeling of distributed environments.
Another interesting variant is multi-skilled scheduling, in which modes require skills instead of explicit resources, with resources mastering one or multiple skills \citep{snauwaert_classification_2023}.
Second, to improve the performance of our CP solvers, it could be worthwhile to design a matheuristic that implements a metaheuristic (such as large neighborhood search) on top of constraint programming.
While CP solvers already use large neighborhood search under the hood, \citet{kasapidis_unified_2024} show that problem-specific destroy operators could further improve its performance.
Finally, we plan to integrate PyJobshop with other CP solvers.
MiniZinc \citep{stuckey_minizinc_2014} offers a standard modeling interface, which can be used to easily integrate a large number of CP solvers without having to write dedicated implementations as we did for OR-Tools and CP Optimizer.

To conclude, we are excited about promoting the use of constraint programming for scheduling problems and we hope that PyJobShop serves as a valuable tool for both researchers and practitioners in this pursuit.

%% file: appendices/experiments.tex
\section{Alternative CP model for FJSP}\label{app:fjsp_naderi}
We describe here the OR-Tools CP model implemented by \citet{naderi_mixed-integer_2023} for the FJSP.
We use the same notation as in our main paper.
In the FJSP, each job $j \in J$ has a set of tasks $T_j$ which must be processed in sequence.
Let $C^\text{EndBeforeStart}$ define all such timing constraints between every pair of consecutive tasks $i$ and $j$.
The goal is to minimize the makespan.
The CP model is described as follows.
\begin{subequations}
\label{model:fjsp_naderi}
\begin{align}
    \label{cons:fjsp_naderi_objective}
    \min \quad &\max_{t \in T, m \in M_t} \mu_{m}^{\supend}  \\
    \label{cons:fjsp_naderi_select_one_mode_variable}
    &\sum_{m \in M_t} \mu_{m}^\present = 1 & \quad \forall t \in T \\ 
    \label{cons:fjsp_naderi_no_overlap}
    &\NoOverlap(\{\mu_{m}: m \in M^R_r \}) &\quad \forall r \in R \\
    &\mu_{m_i}^\supend \le \mu_{m_k}^\start & \quad \forall (i, k, l) \in C^\text{EndBeforeStart}, m_i \in M_i, m_k \in M_k \label{cons:fjsp_naderi_precedence}
\end{align}
\end{subequations}
Expression~\eqref{cons:fjsp_naderi_objective} minimizes the makespan objective.
Constraints~\eqref{cons:fjsp_naderi_select_one_mode_variable} ensure that for each task, exactly one mode variable is selected.
Constraints~\eqref{cons:fjsp_naderi_no_overlap} ensure that there is no overlap on machines.
Constraints~\eqref{cons:fjsp_naderi_precedence} ensure that the timing constraints between consecutive tasks are respected.
In particular, Constraints~\eqref{cons:fjsp_naderi_precedence} are inefficient, as this defines end-before-start constraints between every pair of modes belonging to the corresponding tasks. 
This can be more concisely expressed using the task interval representation as described in Section~\ref{subsec:model}, possibly resulting in better constraint propagation.